\documentclass[reqno,12pt]{amsart}

\usepackage{amsmath, epsfig, cite}
\usepackage{amssymb}
\usepackage{amsfonts}
\usepackage{latexsym}
\usepackage{amsthm}

\newtheorem{theorem}{Theorem}

\newtheorem{corollary}[theorem]{Corollary}
\newtheorem{conjecture}{Conjecture}
\newtheorem{lemma}{Lemma}

\theoremstyle{remark}

\numberwithin{equation}{section}

\allowdisplaybreaks

\author{Yudong Liu and Xiaoxia Wang$^*$}
\address{Department of Applied Mathematics\\
    Shanghai University \\
    Shanghai 200444, P.\:R.\:China}
\email{lydshdx@163.com (Y. Liu), xiaoxiawang@shu.edu.cn (X. Wang)}
\thanks{This work is supported by National Natural Science Foundations of China (11661032).\\
$\quad^*$ Corresponding author.}

\title[Some $q$-supercongruences from Rahman's summation formula]{Some $q$-supercongruences from Rahman's summation}

\subjclass[2010]{Primary 33D15; Secondary 11A07, 11B65}

\keywords{Basic hypergeometric series; supercongruences; $q$-congruences;
cyclotomic polynomial; Rahman's summation formula;}

\begin{document}

\begin{abstract}
Inspired by the recent work on $q$-congruences and the quadratic summation formula of Rahman, we provide some new $q$-supercongruences.
By taking $q\to 1$ in one of our results, we obtain a new Ramanujan-type supercongruence,
which has the same right-hand side as Van Hamme's (G.2) supercongruence for $p\equiv 1 \pmod 4$.
We also formulate some related challenging conjectures on supercongruences and $q$-supercongruences.
\end{abstract}

\maketitle

\section{Introduction}
In 1913, Ramanujan announced the following infinite sum in his first letter to Hardy
(cf.\ \cite[p.~25, Eq.~(2)]{BR}),
\begin{align}
\sum_{k=0}^\infty(8k+1)\frac{(\frac{1}{4})_k^4}{k!^4}
=\frac{2\sqrt{2}}{\sqrt{\pi}\,\Gamma(\frac 34)^2}, \label{eq:ram1}
\end{align}
where $(a)_n=a(a+1)\cdots(a+n-1)$ denotes the Pochhammer symbol and $\Gamma(x)$ is the Gamma function.
Later, Hardy \cite[p.~495]{Ha}  gave a proof of  this infinite summation. Recently, Chen and Chu\cite[Eq. (4) ]{Chu1}  reproved it by giving it's $q$-analogue as follows:
\begin{equation}
\sum_{k=0}^{\infty}[8k+1] \frac{(q ; q^{4})_{k}^4}{(q^4 ; q^{4})_k^4}q^{2k}
=\frac{(q^5,q^3,q^3,q^3;q^4)_\infty}{(q^2,q^4,q^4,q^4;q^4)_\infty}.\label{eq:Chu1}
\end{equation}
They also proposed a similar $q$-infinite summation\cite[Example(8)]{Chu1}:
\begin{align}
\sum_{k=0}^{\infty} [6k+1] \frac{(q ; q^{4})_{k}(q ; q^{2})_{k}^{3}}{(q^{2} ; q^{2})_{k}(q^{4} ; q^{4})_{k}^{3}} q^{k^{2}+k}=\frac{(q^5,q^3,q^3,q^3;q^4)_\infty}{(q^2,q^4,q^4,q^4;q^4)_\infty},\label{eq;Chu2}
\end{align}
whose right-hand side is the same as \eqref{eq:Chu1}.
According to {\cite[p.~353, Eq.~(I.36)]{GM}},
\begin{align}
\lim _{q\to  1^{-}} \frac{(q ; q)_{\infty}}{(q^{x} ; q)_{\infty}}(1-q)^{1-x}=\Gamma(x). \label{q_gamma}
\end{align}
Letting $q\to  1^-$ in \eqref{eq;Chu2}, Chen and Chu obtained the following Ramanujan-type $\pi$-formula:
\begin{align}
\sum_{k=0}^\infty(6k+1)\frac{(\frac{1}{2})_k^3(\frac{1}{4})_k}{4^k k!^4}
=\frac{2\sqrt{2}}{\sqrt{\pi}\,\Gamma(\frac 34)^2}. \label{eq;Chu3}
\end{align}
In 1997, Van Hamme \cite{Hamme} conjectured that \eqref{eq:ram1} has a nice $p$-adic
analogue:
\begin{align}
\sum_{k=0}^{(p-1)/4}(8k+1)\frac{(\frac{1}{4})_k^4}{k!^4}
\equiv p\frac{\Gamma_p(\frac 12)\Gamma_p(\frac 14)}{\Gamma_p(\frac 34)}
\pmod{p^3},\quad\text{if $p\equiv 1\pmod{4}$},
\label{eq:hamme}
\end{align}
where $p$ is an odd prime and $\Gamma_p(x)$ is the $p$-adic Gamma
function~\cite{Mor}. In 2015, Swisher proved this supercongruence \eqref{eq:hamme}. Not long ago, the authors\cite{LW} reproved it by giving the following $q$-analogue:
\begin{align}
\sum_{k=0}^{m}[8k+1] \frac{(q ; q^{4})_{k}^4}{(q^4 ; q^{4})_k^4}q^{2k} \equiv \frac{(q^{2} ; q^{4})_{(n-1)/{4}}}{(q^{4} ; q^{4})_{(n-1)/4}}[n] q^{(1-n) / 4}  \pmod {[n]\Phi_{n}(q)^{2}}\label{eq;qG2},
\end{align}
where $m=(n-1)/4$ or $n-1$.
Recently, Guo and Schlosser  \cite[Theorems 2]{GS2} proved that, for any even integer $d \ge 4$ and positive integer $n$ with $n\equiv-1 \pmod d$,
\begin{align}
\sum_{k=0}^{n-1}[2dk+1]\frac{(q;q^d)_k^d}{(q^d;q^d)_k^d}q^{d(d-3)k/2}
\equiv 0\pmod{\Phi_n(q)^2}. \label{eq:Guo1}
\end{align}
Setting $d=4$ in the above $q$-supercongruence, we have
\begin{align*}
\sum_{k=0}^{n-1}[8k+1]\frac{(q;q^4)_k^4}{(q^4;q^4)_k^4}q^{2k}
\equiv 0\pmod{\Phi_n(q)^2}\quad\text{for $n\equiv 3 \pmod 4$.}
\end{align*}
Here and in what follows, $(a;q)_n=(1-a)(1-aq)\cdots (1-aq^{n-1})$
is the {\em $q$-shifted factorial}. For convenience, we write
$(a_1,a_2,\ldots,a_m;q)_n=(a_1;q)_n (a_2;q)_n\cdots (a_m;q)_n$ for
$n\ge 0 $ or $n=\infty$. Moreover, $[n]:=[n]_q=(1-q^n)/(1-q)$ denotes the {\em
$q$-integer}. Furthermore, $\Phi_n(q)$ is the $n$-th {\em cyclotomic
polynomial} in $q$, which is defined as
\begin{align*}
\Phi_n(q)=\prod_{\substack{1\leqslant k\leqslant n\\ \gcd(n,k)=1}}(q-\zeta^k),
\end{align*}
where $\zeta$ is an $n$-th primitive root of unity.

In the past few years, Ramanujan's and Ramanujan-type $\pi$-formulas
have been generalized to the $q$-world by many authors with various methods (see, for example,
\cite{GL18,G20,GuoZu1,Chu2,Chu1,Wei,Gui,HKS,Sun}).
At the same time, the Ramanujan-type (super)congruences and their $q$-analogues have gradually aroused the research interest of many authors
(see \cite{Gorodetsky,Guo-Dwork,Guo-result,GS2,GS3,GS4,GuoZu,LP,NP,Tauraso2,WY,Zud2009,Zudilin}).
Nowadays, the `creative microscoping' method recently introduced by Guo and Zudilin\cite{GuoZu} makes the research of $q$-supercongruences
much easier than before.

In this paper, we shall give two $q$-supercongruences on truncated forms of \eqref{eq;Chu2}  via `creative microscoping' method as follows.
\begin{theorem}\label{thm:1}
Let $n\equiv 1\pmod {4}$ be a positive integer. Then, modulo $[n]\Phi_{n}(q)^2$,
\begin{align}
&\sum_{k=0}^{(n-1) / 2}[6 k+1] \frac{\left(q ; q^{4}\right)_{k}\left(q ; q^{2}\right)_{k}^{3}}{\left(q^{2} ; q^{2}\right)_{k}\left(q^{4} ; q^{4}\right)_{k}^{3}} q^{k^{2}+k} \notag \\&\equiv \begin{cases}
\frac{\left(q^{2} ; q^{4}\right)_{\frac{(n-1)}{4}}}{\left(q^{4} ; q^{4}\right)_{\frac{(n-1)}{4}}}[n] q^{(1-n) / 4} ,&\quad \text{if $n\equiv 1 \pmod{4}$;} \\0 ,&\quad \text{if $n\equiv 3 \pmod{4}$},
\end{cases}\label{eq;thm1_1}
\end{align}
\begin{align}
&\sum_{k=0}^{n-1}[6 k+1] \frac{\left(q ; q^{4}\right)_{k}\left(q ; q^{2}\right)_{k}^{3}}{\left(q^{2} ; q^{2}\right)_{k}\left(q^{4} ; q^{4}\right)_{k}^{3}} q^{k^{2}+k} \notag \\&\equiv \begin{cases}
\frac{\left(q^{2} ; q^{4}\right)_{\frac{(n-1)}{4}}}{\left(q^{4} ; q^{4}\right)_{\frac{(n-1)}{4}}}[n] q^{(1-n) / 4} ,&\quad \text{if $n\equiv 1 \pmod{4}$;} \\0 ,&\quad \text{if $n\equiv 3 \pmod{4}$}.
\end{cases}\label{eq;thm1_2}
\end{align}
\end{theorem}
Let $n=p \equiv 1\pmod 4$ be a prime and let $q\to  1$ in \eqref{eq;thm1_1}. We have
\begin{equation}
\sum_{k=0}^{(p-1)/2}(6k+1)\frac{(\frac{1}{2})_k^3(\frac{1}{4})_k}{4^k{k!}^4}\equiv \frac{(\frac{1}{2})_{(p-1) / 4}}{(1)_{(p-1) / 4}} p
\pmod {p^3}\label{eq;thm1_3}.
\end{equation}
Notice that, the authors \cite{LW} proved the following congruence holds for prime $p\equiv 1\pmod 4 $,
 \begin{align}
\frac{(\frac{1}{2})_{(p-1) / 4}}{(1)_{(p-1) / 4}} p
\equiv \frac{\Gamma_p(\frac 12)\Gamma_p(\frac 14)}{\Gamma_p(\frac 34)}p \pmod{p^3}. \label{eq;liu}
\end{align}
From \eqref{eq;thm1_1}, we immediately deduce that
\begin{equation}
\sum_{k=0}^{(p-1)/2}(6k+1)\frac{(\frac{1}{2})_k^3(\frac{1}{4})_k}{4^k{k!}^4}
\equiv \begin{cases}
\dfrac{\Gamma_p(\frac 12)\Gamma_p(\frac 14)}{\Gamma_p(\frac 34)}p \pmod{p^3} ,& \text{if $p\equiv 1 \pmod{4}$,}\\[10pt] 0  \pmod{p^3} , & \text{if $p\equiv 3 \pmod{4}$}.
\end{cases}
\label{eq;thm1_4}
\end{equation}
It is a remarkable fact that the infinite Ramanujan-type $\pi$-formulas \eqref{eq:ram1} and \eqref{eq;Chu3} not only
have the same summation, but also have similar $p$-adic analogues. It is easy to see that the
following supercongruences from  \eqref{eq;thm1_2} hold true:
\begin{equation*}
\sum_{k=0}^{p-1}(6k+1)\frac{(\frac{1}{2})_k^3(\frac{1}{4})_k}{4^k{k!}^4}
\equiv \begin{cases}
\dfrac{\Gamma_p(\frac 12)\Gamma_p(\frac 14)}{\Gamma_p(\frac 34)}p \pmod{p^3} ,& \text{if $p\equiv 1 \pmod{4}$,}\\[10pt]
0  \pmod{p^3} , & \text{if $p\equiv 3 \pmod{4}$}.
\end{cases}
\end{equation*}

Motivated by \eqref{eq;qG2} and Theorem \ref{thm:1}, we would like to propose the following challenging conjecture.
\begin{conjecture}
Let $n\equiv 1\pmod {4}$ be a positive integer. Then
\begin{align*}
\sum_{k=0}^{(n-1)/2}[6 k+1] \frac{(q ; q^{4})_{k}(q ; q^{2})_{k}^{3}}{(q^{2} ; q^{2})_{k}(q^{4} ; q^{4})_{k}^{3}} q^{k^{2}+k}
&\equiv \sum_{k=0}^{(n-1)/4}[8k+1] \frac{(q ; q^{4})_{k}^4}{(q^4 ; q^{4})_k^4}q^{2k}\pmod{[n] \Phi_{n}(q)^{3}},\\
\sum_{k=0}^{n-1}[6 k+1] \frac{(q ; q^{4})_{k}(q ; q^{2})_{k}^{3}}{(q^{2} ; q^{2})_{k}(q^{4} ; q^{4})_{k}^{3}} q^{k^{2}+k}
&\equiv \sum_{k=0}^{n-1}[8k+1] \frac{(q ; q^{4})_{k}^4}{(q^4 ; q^{4})_k^4}q^{2k}\pmod{[n] \Phi_{n}(q)^{4}}.
\end{align*}
\end{conjecture}

The rest of the paper is organized as follows. We shall prove Theorem \ref{thm:1} in the next section  with the help of Rahman's summation formula.
Generalizations of Theorem 1  will be given in Section 3. Finally, in Section 4, we will
propose  similar  $q$-congruences, a Ramanujan-type $\pi$-series and a related conjecture on $q$-supercongruences.

\section{Proof of Theorem 1 }
We need the following lemmas in our proof of
Theorems \ref{thm:1}.
\begin{lemma}\label{lem:one}
Let $n$, $d$ be  positive integers with $n>1$ and $n\equiv 1 \pmod {2d} $. Then
\begin{equation}
\sum_{k=0}^{(n-1)/d}[3dk+1]\frac{(q;q^{2d})_k(q^{d-1};q^d)_k(aq,q/a;q^d)_k}{(q^d;q^d)_k(q^{d+2};q^{2d})_k(aq^{2d},q^{2d}/a;q^{2d})_k}q^{d(k^{2}+k)/2}
\equiv 0\pmod{\Phi_{n}(q)}\label{eq;lemma1}.
\end{equation}

\end{lemma}
\begin{proof}
We will make use of a quadratic summation formula of Rahman \cite[Eq. (1.8)]{Rahman}, stated as follow:
\begin{equation}
\sum_{k=0}^{\infty} \frac{1-a q^{3 k}}{1-a} \frac{(a ; q^{2})_{k}(q a / b d ; q)_{k}(b, d ; q)_{k}}{(q ; q)_{k}(q b d ; q^{2})_{k}(a q^{2} / b, a q^{2} / d ; q^{2})_{k}} q^{(k^{2}+k)/2}
= \frac{(a q^{2}, q b, q d, a q^{2} / b d ; q^{2})_{\infty}}{(q, q^{2} a / b, q^{2}a/d, qbd ; q^{2})_{\infty}}.\label{eq;Rm}
\end{equation}
By letting $a=b=d=q^{\frac{1}{2}}$ and $q\to  q^{2}$ in \eqref{eq;Rm}, Chen and Chu \cite[Eq. (4) ]{Chu1} proved \eqref{eq:Chu1}.
Here, setting $q\to  q^{d}$, $a=q^{1-n}$,  $b=aq$ and $d=q/a$ in \eqref{eq;Rm}, then for $n\equiv 1\pmod {2d}$, we have
\begin{align}
&\sum_{k=0}^{(n-1)/2d}\frac{1-q^{3dk+1-n}}{1-q^{1-n}}\frac{(q^{1-n};q^{2d})_k(q^{d-1-n};q^d)_k(aq,q/a;q^d)_k}{(q^d;q^d)_k(q^{d+2};q^{2d})_k(aq^{2d-n},q^{2d-n}/a;q^{2d})_k}q^{d(k^2+k)/2} \notag \\
&=\frac{(q^{2d-1-n},q^{2d+1-n},q^{d+1}/a,aq^{d+1};q^{2d})_\infty}{(q^d,q^{d+2},q^{2d-n}/a,aq^{2d-n};q^{2d})_\infty}.\label{eq;lemma2}
\end{align}
It is easy to check that $(q^{2d+1-n};q^{2d})_\infty =0$ for $n\equiv 1 \pmod{2d}$, thus the right-hand side of the above $q$-congruence is equal to $0$.
Since $q^{n} \equiv 1\pmod {\Phi_{n}(q)}$ and for $(n-1)/2d < k \le (n-1)/d$, $(q;q^{2d})_k/((q^d;q^d)_k(q^{d+2};q^{2d})_k) \equiv 0 \pmod{\Phi_{n}(q)}$, we get
\begin{equation*}
\sum_{k=0}^{(n-1)/d}[3dk+1]\frac{(q;q^{2d})_k(q^{d-1};q^d)_k(aq,q/a;q^d)_k}{(q^d;q^d)_k(q^{d+2};q^{2d})_k(aq^{2d},q^{2d}/a;q^{2d})_k}q^{d(k^{2}+k)/2}\equiv 0 \pmod{\Phi_{n}(q)},
\end{equation*}
as desired.
\end{proof}
\begin{lemma}\label{lemma:two}
Let $n$, $d$ be  positive integers with $n>1$ and $n\equiv -1 \pmod {2d} $. Then
\begin{equation}
\sum_{k=0}^{(n+1)/d-1}[3dk+1]\frac{(q;q^{2d})_k(q^{d-1};q^2)_k(aq,q/a;q^2)_k}{(q^d;q^d)_k(q^{d+2};q^{2d})_k(aq^{2d},q^{2d}/a;q^{2d})_k}q^{d(k^{2}+k)/2}\equiv 0\pmod{\Phi_{n}(q)}. \label{eq;lemma3}
\end{equation}
\end{lemma}
\begin{proof}
Take same parameters transformation in \eqref{eq;Rm} as the proof of  Lemma \ref{lem:one}, the left-hand side of \eqref{eq;Rm} is truncated at $(n+1)/d-1$ place and the right-hand side is same to \eqref{eq;lemma2}, which is equal to $0$ for $n\equiv -1 \pmod {2d}$. Then applying $q^{n} \equiv 1\pmod {\Phi_{n}(q)}$, we immediately obtain \eqref{eq;lemma3}. 
\end{proof}
In order to prove Theorem \ref{thm:1}, we also need to establish the following parametric generalization.
\begin{theorem}\label{thm:2}
Let $n$ be a positive odd integer. Then, modulo $\Phi_{n}(q)(1-a q^{n})(a-q^{n})$,
\begin{align}
&\sum_{k=0}^{(n-1) / 2}[6 k+1] \frac{(q ; q^{4})_{k}(q ; q^{2})_{k}(a q, q / a ; q^{2})_{k}}{(q^{2} ; q^{2})_{k}(q^{4} ; q^{4})_{k}(a q^{4}, q^{4} / a ; q^{4})_{k}} q^{k^{2}+k} \notag \\
&\quad\equiv \begin{cases}
\dfrac{(q^{2} ; q^{4})_{(n-1)/4}}{(q^{4} ; q^{4})_{(n-1)/4}}[n] q^{(1-n) / 4} ,& \text{if $n\equiv 1 \pmod{4}$,} \\[10pt]
0 ,& \text{if $n\equiv 3 \pmod{4}$}.
\end{cases}\label{eq;thm2_1}
\end{align}
\end{theorem}
\begin{proof}
For $a=q^{n}$ or $a=q^{-n}$, the left-hand side of \eqref{eq;thm2_1} is equal to
\begin{align*}
&\sum_{k=0}^{(n-1) / 2}[6 k+1] \frac{(q ; q^{4})_{k}(q ; q^{2})_{k}(q^{1-n}, q ^{1+n} ; q^{2})_{k}}{(q^{2} ; q^{2})_{k}(q^{4} ; q^{4})_{k}(q^{4-n}, q^{4+n} ; q^{4})_{k}} q^{k^{2}+k}\notag\\[5pt]
&\quad=\frac{(q^{5}, q^{3}, q^{3-n}, q^{3+n} ; q^{4})_{\infty}}{(q^{2},  q^{4}, q^{4-n} , q^{4+n} ; q^{4})_{\infty}},
\end{align*}
which is  Rahman's summation formula \eqref{eq;Rm} with the parameter substitutions $q \mapsto q^{2}, a=q, b=q^{1-n}$ and $d=q^{1+n}$.
It's easy to see that
\begin{equation*}
\frac{(q^{5}, q^{3}, q^{3-n}, q^{3+n} ; q^{4})_{\infty}}{(q^{2},  q^{4}, q^{4-n} , q^{4+n} ; q^{4})_{\infty}}=
\begin{cases}
\dfrac{(q^{2} ; q^{4})_{(n-1)/4}}{(q^{4} ; q^{4})_{(n-1)/4}}[n] q^{(1-n) / 4}, & \text{if $n\equiv 1 \pmod{4}$,} \\[10pt]
0 ,& \text{if $n\equiv 3 \pmod{4}$}.
\end{cases}
\end{equation*}
This proves the truth of $q$-congruences \eqref{eq;thm2_1} modulo $(1-aq^{n})(a-q^n)$. Finally, since $\Phi_{n}(q)$, $a-q^n$ and $1-aq^n$ are pairwise relatively prime polynomials,
applying Lemmas \ref{lem:one} and \ref{lemma:two} with $d=2$, we finish the proof of Theorem \ref{thm:2}.
\end{proof}

\begin{proof}[Proof of Theorem 1.] For $k$ in the range $0\le k \le (n-1)/2$, the limits of the denominators on both sides of \eqref{eq;thm2_1} related  to $a$
are relatively prime to $\Phi_{n}(q)$ when $a\to 1$. On the other hand, the limit $(1-aq^n)(a-q^n)$ as  $a\to 1$ contains the factor $\Phi_{n}(q)^2$.
Therefore, the limiting case $a\to 1$ of \eqref{eq;thm2_1} leads to \eqref{eq;thm1_1} are true modulo $\Phi_{n}(q)^3$.
Observing that $ (q ; q^{4})_{k}(q ; q^{2})_{k}^{3}/({(q^{2} ; q^{2})_{k}(q^{4} ; q^{4})_{k}^{3}})\equiv 0 \pmod {\Phi_{n}(q)^{3}}$  for $k$ in the range
$(n-1)/2 \le k \le (n-1)$, we immediately get
\begin{align}
&\sum_{k=0}^{n-1}[6 k+1] \frac{(q ; q^{4})_{k}(q ; q^{2})_{k}^{3}}{(q^{2} ; q^{2})_{k}(q^{4} ; q^{4})_{k}^{3}} q^{k^{2}+k} \notag \\
&\quad\equiv \begin{cases}
\dfrac{(q^{2} ; q^{4})_{(n-1)/4}}{(q^{4} ; q^{4})_{(n-1)/4}}[n] q^{(1-n) / 4} ,& \text{if $n\equiv 1 \pmod{4}$} \\[10pt]
0 ,& \text{if $n\equiv 3 \pmod{4}$}
\end{cases} \pmod{\Phi_{n}(q)^3}.\label{eq;thm2_2}
\end{align}

It remains to prove that the following  congruences also hold modulo $[n]$, i.e.,
\begin{align}
&\sum_{k=0}^{m}[6 k+1] \frac{(q ; q^{4})_{k}(q ; q^{2})_{k}^{3}}{(q^{2} ; q^{2})_{k}(q^{4} ; q^{4})_{k}^{3}} q^{k^{2}+k} \notag \\
&\quad\equiv \begin{cases}
\dfrac{(q^{2} ; q^{4})_{(n-1)/4}}{(q^{4} ; q^{4})_{(n-1)/4}}[n] q^{(1-n) / 4} ,& \text{if $n\equiv 1 \pmod{4}$,} \\[10pt]
0 ,& \text{if $n\equiv 3 \pmod{4}$}.
\end{cases}\pmod{[n]},\label{eq;thm2_4}
\end{align}
where $m$ equals $(n-1)/2$ or $n-1$.
For $n>1$, let $\zeta \not=1 $ be an $n$-th unity root, not necessarily primitive. Then $\zeta$ must be a primitive  $m_1$-th root of unity with $m_1 | n$.
Let $c_q(k)$ denote the $k$-th term on the left-hand side in \eqref{eq;thm2_1}, i.e,
\begin{equation}
c_q(k)=[6k+1]\frac{(q;q^4)_k(q;q^2)_k(aq,q/a;q^2)_k}{(q^2;q^2)_k(q^4;q^4)_k(aq^4,q^4/a;q^4)_k}q^{k^{2}+k}.
\end{equation}
 Applying Lemmas \ref{lem:one} and \ref{lemma:two} with $d=2$, $n=m_1$, and the fact that $(q;q^{2})_k/(q^2;q^2)_k \equiv 0 \pmod{\Phi_{m_1}(q)}$ for $(m_1-1)/2 < k \le m_1-1$, we have
\begin{align*}
\sum_{k=0}^{(m_{1}-1)/2} c_{\zeta}(k)=\sum_{k=0}^{m_{1}-1} c_{\zeta}(k)=0.
\end{align*}
For $0\le k\le m_1-1$, the following limit can be calculated:
\begin{equation*}
\lim _{q\to  \zeta} \frac{c_{q}(l m_{1}+k)}{c_{q}(l m_{1})}=c_{\zeta}(k).
\end{equation*}
Thus, we get
\begin{align}
\sum_{k=0}^{(n-1) / 2} c_{\zeta}(k)=\sum_{\ell=0}^{(n / m_1-3) / 2} c_{\zeta}(\ell m_1) \sum_{k=0}^{m_1-1} c_{\zeta}(k)+\sum_{k=0}^{(m_1-1) / 2} c_{\zeta}((n-m_1) / 2+k)=0,
\end{align}
and
\begin{align}
\sum_{k=0}^{n-1} c_{\zeta}(k)=\sum_{\ell=0}^{n / m_1-1} \sum_{k=0}^{m_1-1} c_{\zeta}(\ell m_1+k)=\sum_{\ell=0}^{n / m_1-1} c_{\zeta}(\ell m_1) \sum_{k=0}^{m_1-1} c_{\zeta}(k)=0.
\end{align}
Noting that
\[\prod_{m_1 | n, m_1>1} \Phi_{m_1}(q)=[n],\]
we establish \eqref{eq;thm2_4}. Since the least common multiple of $[n]$ and $\Phi_{n}(q)^3$ is $[n]\Phi_{n}(q)^2$, we complete the proof of the theorem.
\end{proof}


\section{Generalizations of Theorem 1 }\label{sec:final}
In this section, we first shall give a generalization of $q$-supercongruences \eqref{eq;thm1_1} for $n\equiv 1 \pmod {4}$ modulo $\Phi_{n}(q)^3$.
\begin{theorem}\label{thm:3}
Let $d,n>1$ be integers with $ n\equiv 1 \pmod {2d}$. Then
\begin{align}
&\sum_{k=0}^{(n-1)/d}[3 d k+1] \frac{(q ; q^{2 d})_{k}(q^{d-1} ; q^{d})_{k}( q  ; q^{d})_{k}^2}{(q^{d} ; q^{d})_{k}(q^{d+2} ; q^{2 d})_{k}(q^{2 d}  ; q^{2d} )_{k}^2} q^{d({k^{2}+k})/{2}}\notag\\
&\quad\equiv \frac{(q^{d} ; q^{2 d})_{(n-1)/(2d)}}{(q^{d+2} ; q^{2 d})_{(n-1)/(2d)}}[n] q^{-(n-1)(d-1) / 2 d}\pmod {\Phi_{n}(q)^3}, \label{eq;thm3_1}
\end{align}
\begin{align}
&\sum_{k=0}^{n-1}[3 d k+1] \frac{(q ; q^{2 d})_{k}(q^{d-1} ; q^{d})_{k}( q  ; q^{d})_{k}^2}{(q^{d} ; q^{d})_{k}(q^{d+2} ; q^{2 d})_{k}(q^{2 d}  ; q^{2d} )_{k}^2} q^{d({k^{2}+k})/{2}}\notag\\
&\quad\equiv \frac{(q^{d} ; q^{2 d})_{(n-1)/(2d)}}{(q^{d+2} ; q^{2 d})_{(n-1)/(2d)}}[n] q^{-(n-1)(d-1) / 2 d}\pmod {\Phi_{n}(q)^3}. \label{eq;thm3_2}
\end{align}
\end{theorem}

Similarly, in order to prove Theorem \ref{thm:3}, we first establish the following parametric generalization of \eqref{eq;thm3_1}.
\begin{theorem}\label{thm:4}
Let $d,n>1$ be positive integers with $ n\equiv 1 \pmod {2d}$. Then, modulo $\Phi_{n}(q)(1-aq^n)(a-q^n) $

\begin{align}
&\sum_{k=0}^{(n-1) / d}[3 d k+1] \frac{(q ; q^{2 d})_{k}(q^{d-1} ; q^{d})_{k}(a q, q / a ; q^{d})_{k}}{(q^{d} ; q^{d})_{k}(q^{d+2} ; q^{2 d})_{k}(a q^{2 d}, q^{2 d}/a; q^{2d} )_{k}} q^{d({k^{2}+k})/{2}}\notag \\
&\quad\equiv \frac{(q^{d} ; q^{2 d})_{(n-1)/(2d)}}{(q^{d+2} ; q^{2 d})_{(n-1)/(2d)}}[n] q^{-(n-1)(d-1) / 2 d}. \label{eq;thm4_1}
\end{align}
\end{theorem}

\begin{proof}
For $a=q^{n}$ or $a=q^{-n}$, the left-hand side of \eqref{eq;thm4_1} is equal to
\begin{equation}
\sum_{k=0}^{(n-1) / d}[3 d k+1] \frac{(q ; q^{2 d})_{k}(q^{d-1} ; q^{d})_{k}(q^{1+n}, q^{1-n} ; q^{d})_{k}}{(q^{d} ; q^{d})_{k}(q^{d+2} ; q^{2 d})_{k}(q^{2d+n}, q^{2d-n} ; q^{2d} )_{k}} q^{d({k^{2}+k})/{2}},
\end{equation}
which is the left-hand side of \eqref{eq;Rm} with the parameters replaced by $q \mapsto q^{d}, a=q, b=q^{1-n}$ and $d=q^{1+n}$. Then, the right-hand side of \eqref{eq;Rm} can be written as
\begin{equation*}
 \frac{(q^{2d-1},q^{2d+1},q^{d+1-n},q^{d+1+n};q^{2d})_\infty}{(q^d,q^{d+2},q^{2d-n},q^{2d+n};q^{2d})_\infty}
 =\frac{(q^{d} ; q^{2 d})_{(n-1)/(2d)}}{(q^{d+2} ; q^{2 d})_{(n-1)/(2d)}}[n] q^{-(n-1)(d-1) / 2 d} .
\end{equation*}
This proves that the $q$-congruence \eqref{eq;thm4_1} holds  modulo $1-aq^n$ and $a-q^n$. Moreover, applying Lemma \ref{lem:one}, we can prove that \eqref{eq;thm4_1} is also true modulo $\Phi_{n}(q)$.
\end{proof}
\begin{proof}[Proof of Theorem \ref{thm:3}.]
Similarly as before, letting $a\to 1$,  we can see the denominator of  \eqref{eq;thm4_1} is relatively to $\Phi_{n}(q)$. Also, the limit of $(1-aq^n)(a-q^n)$ as $a\to 1$ has factor $\Phi_{n}(q)^2$. This means that \eqref{eq;thm3_1} modulo $\Phi_{n}(q)^3$ holds true. On the other hand, since $ (q ; q^{2d})_{k}(q ; q^{d})_{k}^{2}/((q^{d} ; q^{d})_{k}(q^{2d} ; q^{2d})_{k}^{2})\equiv 0 \pmod {\Phi_{n}(q)^{3}}$  for $k$ in the range
$(n-1)/d \le k \le (n-1)$, the $q$-supercongruence \eqref{eq;thm3_2} directly follows from \eqref{eq;thm3_1}.
\end{proof}

We also have the following generalization of \eqref{eq;thm1_1} modulo $\Phi_n(q)^2$ for $n\equiv 3 \pmod 4$.
\begin{theorem}\label{thm:5}
Let $d$ and $n>1$ be positive integers  with $ n\equiv d+1 \pmod {2d}$. Then
\begin{align}
\sum_{k=0}^{(n-1)/d}[3 d k+1] \frac{(q ; q^{2 d})_{k}(q^{d-1} ; q^{d})_{k}( q  ; q^{d})_{k}^2}{(q^{d} ; q^{d})_{k}(q^{d+2} ; q^{2 d})_{k}(q^{2 d}  ; q^{2d} )_{k}^2} q^{d({k^{2}+k})/{2}}\equiv 0 \pmod{\Phi_{n}(q)^2}  \label{eq;thm5_1}
\end{align}
\end{theorem}
\begin{proof}
The proof is the same as that of  Theorem \ref{thm:4} and we need to establish the following one more parameter extension:
\begin{align}
\sum_{k=0}^{(n-1) / d}[3 d k+1] \frac{(q ; q^{2 d})_{k}(q^{d-1} ; q^{d})_{k}(a q, q / a ; q^{d})_{k}}{(q^{d} ; q^{d})_{k}(q^{d+2} ; q^{2 d})_{k}(a q^{2 d}, q^{2 d}/a; q^{2d} )_{k}} q^{d({k^{2}+k})/{2}} \equiv 0 \notag \\
\pmod {(a-q^n)(1-aq^n)} \label{eq;thm5_2}.
\end{align}
Noticing that, for  $ n\equiv d+1 \pmod {2d}$,
\begin{align*}
\frac{(q^{2d+1},q^{2d-1},q^{d+1-n},q^{d+1+n};q^{2d})_\infty}{(q^d,q^{d+2},q^{2d-n},q^{2d+n};q^{2d})_\infty}=0,
\end{align*}
we immediately obtain \eqref{eq;thm5_2} by taking  $q \mapsto q^{d}, a=q, b=q^{1-n}$ and $d=q^{1+n}$ in Rahman's summation \eqref{eq;Rm}. Then letting $a\to 1$, we are led to \eqref{eq;thm5_1}.

\end{proof}
Letting $n = p$ be an odd prime, $d = 4$ and $q\to  1$ in Theorem \ref{thm:5}, we get the following result.
\begin{corollary}
Let $p \equiv 5 \pmod 8                                                                                                                                                                                                                                                                                                                                                                                                                                                                                                                                                                                                                                                                                                                                                                                                                                                                                                                                                                                                                                                                                                                                                                                                                                                                                                                                                                                                                                                                                                                                                                                                                                                                                                                                                                                                                                                                                                                                                                                                                                                                                                                                                                                                                                                                                                                                                                                                                                                                                                                                                                                                                                                                                                                                                                                                                                                                                                                                                                                                                                                                                                                                                                                                                                                                                                                                                                                                                                                                                                                                                                                                                                                                                                                                                                                                                                                                                                                                                                                                                                                                                                                                                                                                                                                                                                                                                                                                                                                                                                                                                                                                                                                                                                                                                                                                                                                                                                                                                                                                                                                                                                                                                                                                                                                                                                                                                                                                                                                                                                                                                                                                                                                                                                                                                                                                                                                                                                                                                                                                                                                                                                                                                                                                                                                                                                                                                                                                                                                                                                                                                                                                                                                                                                                                                                                                                                                                                                                                                                                                                                                                                                                                                                                                                                                                                                                                                                                                                                                                                                                                                                                                                                                                                                                                                                                                                                                                                                                                                                                                                                                                                                                                                                                                                                                                                                                                                                                                                                                                                                                                                                                                                                                                                                                                                                                                                                                                                                                                                                                                                                                                                                                                                                                                                                                                                                                                                                                                                                                                                                                                                                                                                                                                                                                                                                                                                                                                                                                                                                                                                                                                                                                                                                                                                                                                                                                                                                                                                                                                                                                                                                                                                                                                                                                                                                                                                                                                                                                   $ be a prime. Then
\begin{align}
\sum_{k=0}^{(p-1)/4}(12 k+1) \frac{(\frac{1}{4})_{k}^{2}(\frac{1}{8})_{k}}{4^{k} k !^{3}} \equiv 0 \pmod {p^2}.
\end{align}
\end{corollary}
By the numerical calculations, we believe that the following stronger version is also true.
\begin{conjecture}
Let $p \equiv 5 \pmod 8 $ be a prime. Then
\begin{align*}
 \sum_{k=0}^{p-1}(12 k+1) \frac{(\frac{1}{4})_{k}^{2}(\frac{1}{8})_{k}}{4^{k} k !^{3}} \equiv 0 \pmod {p^3}.
\end{align*}
\end{conjecture}
\section{Other Conclusions from Rahman's Summation}\label{sec:final}
In this section, we present similar $q$-supercongruences and a $\pi$-formula from Rahman's quadratic summation formula \eqref{eq;Rm}.
\begin{theorem}\label{thm:7}
Let $n$ and $d \not= 1$ be positive  integers. Then, module $\Phi_{n}(q)^2$,
\begin{align}
&\sum_{k=0}^{(n-1)/d}[3dk-1]\frac{(q^{-1};q^{2d})_k(q^{d-3};q^d)_k(q;q^d)_k^2}{(q^d;q^d)_k(q^{d+2};q^{2d})_k(q^{2d-2};q^{2d})_k^2}q^{d(k^2+k)/2+1} \notag\\[5pt]
&\quad\equiv  \begin{cases}
- \dfrac{(q^{d} ; q^{2d})_{(n-1)/(2d)}(q^{2d-1} ; q^{2d})_{(n-1)/(2d)}}{(q^{3} ; q^{2d})_{(n-1)/(2d)}(q^{d+2} ; q^{2d})_{(n-1)/(2d)}} q^{(3-d)(n-1)/2d}, & \text{if $n\equiv 1 \pmod{2d}$,} \\[10pt]
0, & \text{if $n\equiv d+1 \pmod{2d}$}.
\end{cases}
\end{align}

\end{theorem}\

\begin{proof}
We first establish the following congruence module $(a-q^n)(1-aq^n)$:
\begin{align}
&\sum_{k=0}^{(n-1)/d}[3dk-1]\frac{(q^{-1};q^{2d})_k(q^{d-3};q^d)_k(q/a;q^d)_k(aq;q^d)_k}{(q^d;q^d)_k(q^{d+2};q^{2d})_k(aq^{2d-2};q^{2d})_k(q^{2d-2}/a;q^{2d})_k}q^{d(k^2+k)/2+1}\notag \\[5pt]
&\quad\equiv  \begin{cases}
- \dfrac{(q^{d} ; q^{2d})_{(n-1)/(2d)}(q^{2d-1} ; q^{2d})_{(n-1)/(2d)}}{(q^{3} ; q^{2d})_{(n-1)/(2d)}(q^{d+2} ; q^{2d})_{(n-1)/(2d)}} q^{(3-d)(n-1)/2d} , & \text{if $n\equiv 1 \pmod{2d}$,} \\[10pt]
0 , & \text{if $n\equiv d+1 \pmod{2d}$}.
\end{cases}\label{eq;thm7_pf}
\end{align}
Letting $q \mapsto q^{d}, a=q^{-1}, b=q^{1-n}$ and $d=q^{1+n}$ in \eqref{eq;Rm}, we obtain
\begin{align*}
&\sum_{k=0}^{(n-1)/d}[3dk-1]\frac{(q^{-1};q^{2d})_k(q^{d-3};q^d)_k(q^{1-n};q^d)_k(q^{1+n};q^d)_k}{(q^d;q^d)_k(q^{d+2};q^{2d})_k(q^{2d-2+n};q^{2d})_k(q^{2d-2-n};q^{2d})_k}q^{d(k^2+k)/2+1}\\[5pt]
&\quad=-\frac{(q^{2d-1},q^{2d-3},q^{d+1-n},q^{d+1+n};q^{2d})_\infty}{(q^d,q^{d+2},q^{2d-2-n},q^{2d-2+n};q^{2d})_\infty}\\[5pt]
&\quad=
\begin{cases}
- \dfrac{(q^{d} ; q^{2d})_{(n-1)/(2d)}(q^{2d-1} ; q^{2d})_{(n-1)/(2d)}}{(q^{3} ; q^{2d})_{(n-1)/(2d)}(q^{d+2} ; q^{2d})_{(n-1)/(2d)}} q^{(3-d)(n-1)/2d}, & \text{if $n\equiv 1 \pmod{2d}$,} \\[10pt]
0, & \text{if $n\equiv d+1 \pmod{2d}$.}
\end{cases}
\end{align*}
This proves that the $q$-congruences generalization \eqref{eq;thm7_pf} with one more parameter hold true. The proof then follows by taking  limits as $a\to 1$ in \eqref{eq;thm7_pf}.
\end{proof}
Specially, when $d=2$ in Theorem \ref{thm:7}, modulo $\Phi_{n}(q)^2$, we have
\begin{align}
&\sum_{k=0}^{(n-1)/2}[6k-1]\frac{(q^{-1};q^{4})_k(q^{-1};q^2)_k(q;q^2)_k^2}{(q^2;q^2)_k(q^{4};q^{4})_k(q^{2};q^{4})_k^2}q^{k^2+k+1} \notag \\
&\quad\equiv \begin{cases}
-\dfrac{(q^{2} ; q^{4})_{(n-1)/4}}{(q^{4} ; q^{4})_{(n-1)/4}} q^{(n-1) / 4}, & \text{if $n\equiv 1 \pmod{4}$,} \\[10pt]
0, & \text{if $n\equiv 3 \pmod{4}$}.
\end{cases}
\end{align}
Furthermore, we believe that the following stronger $q$-supercongruences are true.
\begin{conjecture}
Let $n$ be a positive odd integer. Then, modulo $[n]\Phi_{n}(q)^2$,
\begin{align}
&\sum_{k=0}^{(n-1)/2}[6k-1]\frac{(q^{-1};q^{4})_k(q^{-1};q^2)_k(q;q^2)_k^2}{(q^2;q^2)_k(q^{4};q^{4})_k(q^{2};q^{4})_k^2}q^{k^2+k+1} \notag \\
&\quad\equiv \begin{cases}
-\dfrac{(q^{2} ; q^{4})_{(n-1)/4}}{(q^{4} ; q^{4})_{(n-1)/4}} q^{(n-1) / 4} ,& \text{if $n\equiv 1 \pmod{4}$,} \\[10pt]
0,,& \text{if $n\equiv 3 \pmod{4}$}.
\end{cases}
\end{align}
\end{conjecture}
Letting $ n=p$ be a prime, $q\to  1$ in the above $q$-supercongruences,  and applying \eqref{eq;liu} again, we have
\begin{align}
\sum_{k=0}^{(p-1)/2}(6 k-1) \frac{(-\frac{1}{2})_{k}(-\frac{1}{4})_{k}}{4^{k} k !^2}
\equiv \begin{cases}
 -\dfrac{\Gamma_p(\frac 12)\Gamma_p(\frac 14)}{\Gamma_p(\frac 34)}\pmod{p^3} ,
 & \text{if $p\equiv 1 \pmod{4}$,}\\[5pt] 0\qquad \qquad \pmod{p^3},& \text{if $p\equiv 3 \pmod{4}$}.\label{eq;thm7_1}
\end{cases}
\end{align}
Meanwhile, replacing $q$ by $q^2$, $ a=q^{-1}$ and $b=d=q$ in \eqref{eq;Rm}, we get an infinity summation formula:
\begin{align}
\sum_{k=0}^{\infty}[6k-1]\frac{(q^{-1};q^{4})_k(q^{-1};q^2)_k(q;q^2)_k^2}{(q^2;q^2)_k(q^{4};q^{4})_k(q^{2};q^{4})_k^2}q^{k^2+k+1}=-\frac{(q^3,q^3,q^3,q;q^4)_\infty}{(q^2,q^2,q^2,q^4;q^4)_\infty}.
\label{eq;thm7_2}
\end{align}
Applying \eqref{q_gamma} and noting that $\Gamma(x)  \Gamma(1-x)=\pi/\sin (\pi x)$, then letting $q\to 1^{-}$ in \eqref{eq;thm7_2}, we shall obtain the following Ramanujan-type summation of $\pi$ as
\begin{align}
\sum_{k=0}^{\infty}(6 k-1) \frac{(-\frac{1}{2})_{k}(-\frac{1}{4})_{k}}{4^{k} k !^2} &=-\frac{\Gamma(\frac{1}{2})^3\Gamma(1)}{\Gamma(\frac{3}{4})^3\Gamma(\frac{1}{4})}\notag \\
&=-\frac{\sqrt{2\pi}}{2\Gamma^2(\frac{3}{4})}.\label{eq;thm7_3}
\end{align}
Thus, the supercongruence \eqref{eq;thm7_1} may be deemed a nice $p$-adic analogue of \eqref{eq;thm7_3}.

\end{document}